\documentclass[11pt]{article}

\usepackage{amssymb}
\usepackage{amsthm}
\usepackage{amsmath,amsfonts,amssymb}

\allowdisplaybreaks
\usepackage[mathscr]{eucal}
\usepackage{exscale}
\usepackage[square,sort,comma,numbers]{natbib}
\setcitestyle{authoryear,round}

\usepackage{bm}
\usepackage{eqlist}
\usepackage[dvipsnames]{color}
\usepackage[Lenny]{fncychap}
\usepackage{multirow}
\usepackage{graphicx}
\usepackage{algpseudocode}
\usepackage{algorithm}
\usepackage{caption}
\usepackage{hyperref}
\usepackage{framed}
\usepackage{color}
\usepackage{booktabs}
\usepackage{makecell}
\usepackage{framed}
\usepackage{algpseudocode}
\usepackage{indentfirst} 
\usepackage{longtable}

\usepackage{tabu}
\usepackage{tikz}

\usepackage{chngpage}

\usepackage[utf8]{inputenc} 
\usepackage[english]{babel} 
\usepackage[T1]{fontenc}    
\usepackage{amssymb,amsmath,amsfonts} 
\usepackage{braket} 
\usepackage{caption}
\usepackage{algorithm}

\hypersetup{
	colorlinks=true,
	linkcolor=blue,
	citecolor=blue,
	citebordercolor={1 1 0},
}

\newtheorem{theorem}{Theorem}[section]

\newtheorem{definition}[theorem]{Definition}

\font\bigbf=cmbx10 scaled \magstep3

\begin{document}

\title{\bigbf Route planning of mobile sensing fleets for repeatable visits}

\author{Wen Ji$^a$
\quad 
Ke Han$^{b,}\thanks{Corresponding author, e-mail: kehan@swjtu.edu.cn;}$
\quad
Qian Ge$^a$
\\\\
\textit{\small $^a$School of Transportation and Logistics, Southwest Jiaotong University, Chengdu, China}
\\
\textit{\small $^b$School of Economics and Management, Southwest Jiaotong University, Chengdu, China}
}

\maketitle

\begin{abstract}
Vehicle-based mobile sensing is an emerging data collection paradigm that leverages vehicle mobilities to scan a city at low costs. Certain urban sensing scenarios require dedicated vehicles for highly targeted monitoring, such as volatile organic compounds (VOCs, a type of air pollutant) sensing, road surface  monitoring, and accident site investigation. A hallmark of these scenarios is that the points of interest (POIs) need to be repeatedly visited by a set of agents, whose routes should provide sufficient sensing coverage with coordinated overlap at certain important POIs. For these applications, this paper presents the {\it open team orienteering problem with repeatable visits} (OTOP-RV). The adaptive large neighborhood search (ALNS) algorithm is tailored to solve the OTOP-RV considering specific features of the problem. Test results on randomly generated datasets show that: (1) For small cases, the ALNS matches Gurobi in terms of optimality but with shorter computational times; (2) For large cases, the ALNS significantly outperforms the greedy algorithm (by 9.7\% to 25.4\%), and a heuristic based on sequential orienteering problems (by 6\%). Finally, a real-world case study of VOCs sensing is presented, which highlights the unique applicability of the OTOP-RV to such specific sensing tasks, as well as the effectiveness of the proposed algorithms in optimizing the sensing utilities.
\end{abstract}


\section{Introduction}

    Drive-by sensing has gained significant popularity in smart city applications, such as air quality \citep{Song2021, Mahajan2021}, traffic state \citep{ZLMH2014, GQDRJ2022}, noise pollution \citep{Alsina2017}, heat island phenomena \citep{Fekih2021}, and infrastructure health inspection \citep{WBS2014, MCG2022}. Depending on how the fleets are sourced, sensing vehicles can be categorized as crowd-sourced and dedicated vehicles. The former include taxis \citep{oKeeffe2019}, buses \citep{JHL2023b, DH2023}, trams \citep{Saukh2012} and logistics vehicles \citep{dADKKR2020}, and are suitable for large-scale pervasive sensing at relatively low operating costs \citep{JHL2023a}, including air quality, heat island, and traffic state sensing. The latter refer to vehicles that are used for the sole purpose of sensing, with full controllability but at higher costs. Examples of this category include volatile organic compounds (VOCs, a type of air pollutant) monitoring, road surface condition sensing, and accident on-site investigation \citep{Messier2018, GM2020, LXSW2020}. Such application scenarios require detailed route planning for a fleet of dedicated vehicles (DVs).

DV route planning for monitoring tasks with resource constraints can be regarded as TOP \citep{CGW1996b} and its variants proposed in the literature. \cite{GM2020} introduce the concept of mission planning into UAV routing for accident site investigation, and defined a generalized correlated TOP for selecting POIs, taking into account spatial correlations and priorities. \cite{LXSW2020} consider a Min-Time Max-Coverage issue in sweep coverage where a set of UAVs are dispatched to efficiently patrol the POIs in the given area to achieve maximum coverage in minimum time for forest fire monitoring. \cite{DZCM2021} design an actor-centric heterogeneous collaborative reinforcement learning algorithm to schedule different UAVs to maximize sensing coverage, coverage fairness, and cost-effectiveness. \cite{FHWZ2023} design a neural network heuristic for routing multi-UAVs to monitor scattered landslide-prone areas, with mandatory visits on those in poorly stable states. The above research focuses on the route planning problem when there are insufficient agents (DVs or UAVs) to cover all the POIs, so it is necessary to select a subset of POIs based on specific priorities.

It is essential to note that the OP and its variants proposed in the literature require that each POI is visited at most once. In reality, however, there are a significant number of cases where certain POIs need to be visited multiple times, e.g. to collect sufficient information, such as urban mobile sensing (air quality, noise, heat island, etc.), or to perform tasks that require repetition, such as disaster relief (forest fire containment, emergency food supply). Such needs vary among different POIs, and it is important to allow overlap of the agents' routes in a way that meets application-specific requirements.

\begin{figure}[h]
       \centering
       \includegraphics[width=0.9\textwidth]{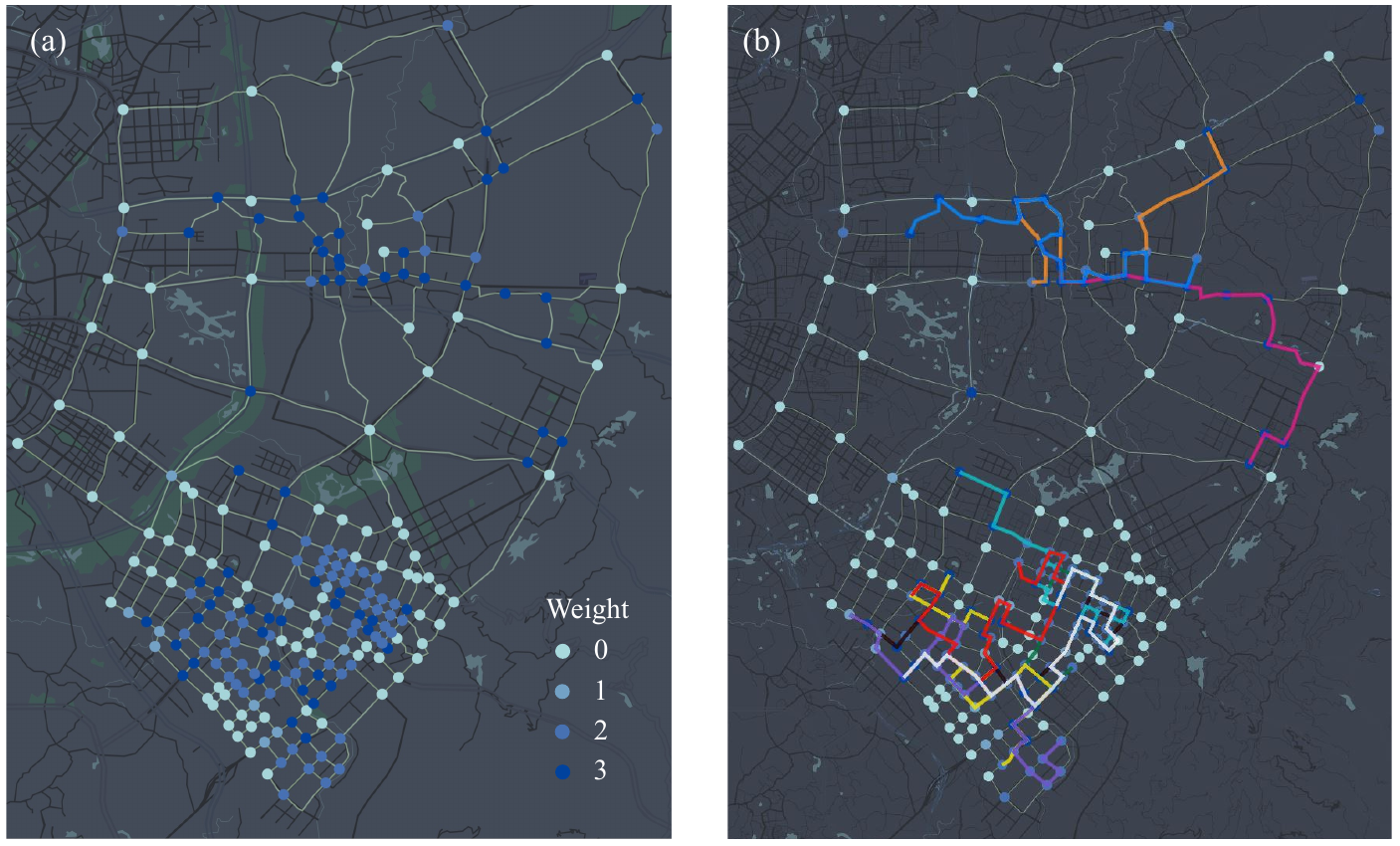}
       \caption{Demonstration of the OTOP-RV solution: (a) A road network where the nodes are POIs with varying weights; (b) the 10 routes generated to visit these POIs.}
       \label{fig_intro}
   \end{figure}

We take, as an example, drive-by sensing of volatile organic compounds (VOCs) performed by a single DV, which needs to consistently monitor a number of factories and collect VOCs measurements. The factories are associated with different sensing importance (weights). The agent needs to execute two routes per day, with a distance budget of 20km per route, which amounts to a total of 10 routes per week (5 workdays). These 10 routes are expected to have a certain overlap since factories with higher sensing weights should be visited more frequently. Figure \ref{fig_intro} shows the 10 routes generated by our proposed algorithm (see Section \ref{subsecRWCS}). No existing variant of TOP is suited for such a situation.
	
    Motivated by this, we propose the {\it open team orienteering problem with repeatable visits} (OTOP-RV), which aims to determine a set of routes, each within a given time/distance budget and covering a subset of POIs, such that the total sensing utility is maximized. Here, `open' means the agent routes do not need to start or end at a fixed depot. Furthermore, if the POIs are arcs instead of vertices, the OTOP-RV can be applied to the route planning of road sprinkler operation \citep{ZSZWZ2022} and road surface condition monitoring \citep{AD2017}, which can be seen as an extension of the classical capacitated arc routing problem (CARP) \citep{GW1981}, with heterogeneous link coverage. 
    
    Specific contributions of this work are as follows: We present a mathematical formulation of the OTOP-RV as a nonlinear integer program, and design an adaptive large neighborhood search (ALNS) algorithm to solve large-scale instances. Test results on randomly generated datasets show that: (1) For small cases, the ALNS algorithm can find the same optimal solutions as Gurobi but in shorter time; (2) For large cases, the ALNS significantly outperform a greedy algorithm (by 9.7\% to 25.4\%) and a heuristic based on sequential OPs (by 6\%). Finally, a real-world case study of mobile VOCs sensing on a road network is presented and formulated as OTOP-RV. Such a case study highlights the unique applicability of OTOP-RV to certain real-world scenarios, as well as the effectiveness of the proposed algorithms in achieving a sensible solution. 
	
	The rest of this article is organized as follows. Section \ref{sec_problem_statement} articulates the OTOP-RV with a mathematical formulation. Sections \ref{sec_ALNS} present the ALNS solution approach to the OTOP-RV. A discussion on model extensions is presented in Section \ref{secDE}. Section \ref{sec_computational_study} presents extensive numerical and application studies. Finally, Section \ref{sec_conclusion} presents some concluding remarks. 

\section{Related work}
	Work related to this study is divided into three parts: DV-based drive-by sensing, team orienteering problem, and multi-visit vehicle routing problem.

\subsection{DV-based drive-by sensing}    
DVs have been used in many drive-by sensing scenarios due to their high controllability and flexibility. The data requirements of drive-by sensing depend on the sensing object and the underlying application. Based on the different sensing requirements, we can categorize DV-based drive-by sensing related research into three groups. The first type, which is typically appropriate for scenarios with low sensing frequency requirements, such as built environment monitoring \citep{LMZLL2023, LL2024}, attempts to achieve complete coverage of the monitored area at the lowest cost. These problems are known as the Traveling Salesman Problem \citep{Bektas2006} or the Chinese postman problem \citep{AR2006}. The second type, which is typically appropriate for scenarios with high response speed, such as search and rescue \citep{YARZM2022}, seeks to maximize the total monitoring reward by selecting a subset of locations to visit with limited monitoring resources. These problems are known as the Orienteering Problem \citep{Golden1987, Gunawan2014} or the Prize-Collecting Vehicle Routing Problem \citep{LSPHZL2019, RHS2021}. The third type is typically appropriate for scenarios with higher sensing frequency requirements, therefore, multiple visits to locations within a period can lead to an increase in sensing reward. The applications include VOC (Volatile Organic Compounds) sensing, road surface condition monitoring, or accident site investigation. However, there is currently no research on vehicle route planning for this type of problem.

This study aims to plan the inspection routes of the DV fleet for the third type of problem, to achieve a balance between monitoring coverage and monitoring frequency. This is a complex problem that requires careful consideration of both the breadth and depth of data collection, ensuring that critical locations are not only covered but also monitored frequently enough to capture relevant changes.

\subsection{Team orienteering problem}

    The orienteering problem (OP) is a routing problem, which determines a subset of nodes (or Points of Interest, POIs) to visit, and the order in which they are visited, so that the total collected profits are maximized within a given time or distance budget \citep{Golden1987, CGW1996a}. The team orienteering problem (TOP) is a natural extension of the OP by generating multiple routes \citep{CGW1996b, DGM2013, KZLC2015}. Typical application scenarios of the TOP include athlete recruiting \citep{CGW1996b}, home fuel delivery \citep{GAD1984}, search and rescue operations \citep{YARZM2022} and tourist trip planning \citep{VV2007}. 

    Several variants of the TOP have been proposed and studied in the literature, including the (team) orienteering problem with time window, in which each POI can only be visited within a given time window \citep{Gambardella2012, LY2012, DLM2015}; the time-dependent orienteering problem, in which the traveling time between two POIs is time-dependent \citep{Verbeeck2014, Gunawan2014}; the stochastic orienteering problem, in which the traveling time and the reward at each POI have stochastic attributes \citep{Ilhan2008, ZOT2014}; the (team) orienteering problem with variable profits, in which the profit depend on the arrival time or the service time at each POI \citep{YFZM2019, YARZM2022, TMT2007}; and the team orienteering problem with mandatory visits, in which some high priority POIs require mandatory visit \citep{PASLL2017, FHWZ2023}. For a more comprehensive review of the OP, its variants, and their applications, we refer the reader to \cite{Vansteenwegen2011} and \cite{Gunawan2016}. 

    In this study, we introduce a new variant of the team orienteering problem, in which each POI may be visited by multiple routes to accumulate rewards, named an open team orienteering problem with repeatable visits.

\subsection{Multi-visit vehicle routing problem}
    Some vehicle routing problems are inspired by real-world applications and involve multiple visits. \cite{HMZ2020} propose a new variant of TOP, called multi-visit TOP with precedence constraints, with an application scenario of on-site assembly of kitchen furniture and appliances. In this problem, each customer has a set of tasks to be performed by a heterogeneous set of technicians according to a predefined order. When a customer is selected, all the tasks have to be performed by possibly different technicians, so the customer may be visited more than once by possibly different technicians. To solve the problem, the authors proposed a compact MILP formulation and a kernel search heuristic algorithm. \cite{PAM2020} study a new workforce scheduling and routing problem where the requested services consist of tasks to be performed over one or more days by teams of workers with different skills. Each customer can be visited more than once, as long as precedence constraints are not violated. To solve the problem, the authors proposed a MILP formulation and an ant colony metaheuristic algorithm. \cite{AISV2023} introduces a vehicle routing problem for a specific application in the monitoring of a water distribution network (WDN). In this problem, multiple technicians must visit a sequence of nodes in the WDN and perform a series of tests to check the quality of water. Some special nodes (i.e., wells) require technicians to first collect a key from a key center and the key must then be returned to the same key center after the test has been performed, thus introducing multiple visits in the routes. To solve the problem, the authors proposed a MILP formulation and an iterated local search heuristic algorithm. In addition, the pickup-and-delivery problem with split loads and transshipments also allows multiple visits to the same locations \citep{BI2017, W2021, WS2021}. This is common when split deliveries are allowed, or multiple pickup and delivery operations can be performed at a single location. 

	In summary, the studies mentioned above all consider multi-visit as a necessary constraint for vehicle routing problems, and the number of visits to different locations has no impact on the rewards. In contrast, this study proposes a new multi-visit vehicle routing problem where the visit reward of each location is closely related to the number of visits.
    
\section{Problem Statement and Model Formulation} \label{sec_problem_statement}

\subsection{Problem description}

	We consider a target monitoring area with a given set of points of interest (POIs), to be visited by a set of agents\footnote{In the following analyses of the orienteering problem and its variants, we use the term agents in replacement of vehicles or vehicle trips.}. To capture the unique characteristic of our problem, which is the repeatable visit, we define a sensing reward function for each POI $i$, which depends on the number of distinct agents that have visited it, denoted $q_i\in\mathbb{Z}_+$\footnote{Visits performed by the same agent can be counted at most once. This stipulation eliminates inefficient routes that intersect themselves.}. Specifically, such a reward function, $\phi(q_i)$, needs to satisfy the following conditions:
\begin{enumerate}
\item $\phi(0)$=0, and $\phi(\cdot)$ is monotonically increasing;
\item The marginal gain of the reward is decreasing: 
$$
\phi(q_{i}+1)-\phi(q_{i})> \phi(q_{j}+1)-\phi(q_{j})\qquad\forall 0\leq q_{i}< q_{j}
$$
\end{enumerate}

\noindent The first condition is straightforward. The second condition implies, as the visits accumulate, the need to visit the POI one more time decreases. It is necessary to avoid over-concentration of visits at a few high-value POIs. The following function, which satisfies these conditions, will be used in this paper.
\begin{equation}\label{eqscfor1}
\phi(q_i)= (q_i)^{\beta}\qquad \beta \in (0, 1)
\end{equation}
\noindent Finally, the objective of the OTOP-RV is to maximize the following weighted sum:
	\begin{equation}
	        \label{eq_sensing_reward}
	        \Phi = \sum_{i \in I}w_{i}\phi(q_i)=\sum_{i \in I}w_{i}(q_{i})^\beta 
	\end{equation}

\noindent where the parameter $w_{i}$ is the weight of POI $i$. Intuitively, when $\beta\to 1$, the agents tend to visit high-weighting POIs more frequently; as $\beta\to 0$, their visits are more evenly distributed among all POIs. The parameters $w_i$'s and $\beta$ should be jointly determined based on the underlying application.

The OTOP-RV is articulated as follows.

    \begin{definition}\textbf{The Open Team Orienteering Problem with Repeatable Visits (OTOP-RV)} Given a set of POIs $I=\{1, ..., n\}$, each with a sensing weight $w_i$ and a reward function $\phi(q_i)=(q_i)^{\beta}$, and a total of $n_K$ agents. The OTOP-RV needs to determine $n_{K}$ routes, each within a fixed time or distance budget, and do not need to start and end at the same depot. The goal is to maximize the total reward collected from visiting these POIs $\sum_{i\in I} w_i (q_i)^{\beta}$.
    \end{definition}

The \underline{open} team orienteering problem allows each route to start and end at any nodes in the network, which is suitable for sensing or surveillance tasks, unlike conventional logistics scenarios. The time or distance budget could be related to the agent's travel limit, or the working hours of operating personnel.

\subsection{Mathematical model}

Table \ref{tab_notations}  lists some key notations used in the model.

\setlength\LTleft{0pt}
\setlength\LTright{0pt}
\begin{longtable}{@{\extracolsep{\fill}}rl}
\caption{Notations and symbols}
\label{tab_notations} 
\\
\hline
\multicolumn{2}{l}{Sets}     
    \\
    \hline
    $I$  & Set of POIs;
    \\
    $K$  & Set of agents;
    \\
    $N$  & Set of vertices in the graph $G$;
    \\
    $A$ & Set of arcs on the graph $G$;
    \\
    \hline
    \multicolumn{2}{l}{Parameters and constants}           
    \\
    \hline
    $n_{K}$  & Number of agents; 
    \\
    $n$  &  Number of POIs;
    \\
    $t_{ij}$  & Travel time from vertex $i$ to vertex $j$;
    \\
    $w_{i}$  & The weight of POI $i$;
    \\
    $\Delta$  & The length of time or distance budget.
    \\
    $\beta$  & The parameter in the reward function \eqref{eqscfor1} that ensures diminishing marginal gain.
    \\
    \hline
    \multicolumn{2}{l}{Auxiliary variables}                
    \\
    \hline
    $y_{i,k}$  & Binary variable that equals $1$ if POI $i$ is visited by agent $k$;
    \\
    $a_{i, k}$  & The visiting time of POI $i$ by agent $k$;
    \\
    $q_{i}$  & The total number of times the POI $i$ is covered by the agents;
    \\
    \hline
    \multicolumn{2}{l}{Decision variables}                 
    \\
    \hline
    $x_{i,j,k}$  & Binary variable that equals $1$ if agent $k$ travels directly from vertex $i$ to vertex $j$.
    \\	
    \hline
\end{longtable}

     We consider a network represented as an undirected graph $G=(N, A)$, where $N$ is the set of vertices, and $A$ is the set of arcs. Although it is not necessary to ensure that the starting and ending points of each route are the same depot in this study, we still need to introduce virtual starting depot $0$ and ending depot $n+1$ to ensure the flow balance and consistency constraints. The arc set consists of two parts: those connecting any pair of POIs, as well as those connecting the virtual starting/ending depots to the POIs: $A = \{(0, i)|i \in I\} \cup \{(i, j)|i,j \in I; t_{ij} \neq \inf \} \cup \{(j, n+1)|j \in I\}$. The full OTOP-RV model is formulated as:

\begin{equation}
        \label{E1}
        \max_{s=\{x_{i,j,k}: \, i,j\in I, k\in K\}} \sum_{i \in I}w_{i}q_{i}^\beta
\end{equation}

    \begin{eqnarray}
        \label{E2}
        \sum_{i \in I}x_{0,i,k}=1  & &  \forall k \in K
        \\
        \label{E3}
        \sum_{i \in I}x_{i,n+1,k}=1  & & \forall k \in K
        \\
        \label{E4}
        \sum_{(i,l) \in A}x_{i,l,k} = \sum_{(l,j) \in A}x_{l,j,k} & & \forall l \in I, \forall k \in K
        \\
        \label{E5}
        \sum_{(i,j) \in A}x_{i,j,k} = y_{i, k}  & & \forall i \in I, \forall k \in K
        \\
        \label{E6}
         q_{i}=\sum_{k \in K}y_{i,k} & & \forall i \in I
         \\
        \label{E7}
        a_{i,k}+t_{ij}-a_{j,k} \leq M(1 - x_{i,j,k}) & & \forall i,j \in N, k \in K
        \\
        \label{E8}
        a_{n+1} \leq \Delta
        \\
        \label{E9}
        x_{i,j,k} \in \{0,1\}  & & \forall i,j \in N, \forall k \in K
        \\
        \label{E10}
        y_{i,k} \in \{0, 1\}  & & \forall i \in I, \forall k \in K
        \\
        \label{E11}
        a_{i, k} \in \mathbb{R}_{+}  & & \forall i \in I, \forall k \in K
    \end{eqnarray}

    The objective function \eqref{E1} maximizes the weighted rewards within time or distance budget $\Delta$. Constraints \eqref{E2} and \eqref{E3} ensure that each agent departs from the virtual starting depot and returns to the virtual ending depot, respectively. Constraint \eqref{E4} ensures flow conservation for each agent. Constraint \eqref{E5} couples variables $x_{i,j,k}$ with indicator variable $y_{i,k}$. Constraint \eqref{E6} calculates the total number of times the POI $i$ is visited by the agents. Constraint \eqref{E7} acts as sub tour elimination constraint \citep{MTZ1960} and ensure that a route cannot visit the same POI multiple times. Constraint \eqref{E8} ensures that the time or distance budget is not exceeded.

\section{ALNS applied to the OTOP-RV} \label{sec_ALNS}
	This section describes a solution framework based on adaptive large neighborhood search (ALNS) for the OTOP-RV. As summarized in \cite{Windras2022}, the ALNS has the following important parts: (1) initial solution $s_0$; (2) a set of destroy operators $\Omega^{-}=\{\Omega_{1}^{-}, ..., \Omega_{|\Omega^{-}|}^{-}\}$; (3) a set of repair operators $\Omega^{+}=\{\Omega_{1}^{+}, ..., \Omega_{|\Omega^{+}|}^{+}\}$; (4) adaptive mechanism; (5) acceptance criterion and (6) termination criterion. In each iteration of ALNS, a single destroy and a single repair operator are selected according to the \textbf{adaptive mechanism}. The destroy operator is first implemented to remove parts of a feasible solution $s$ and the resulting solution is stored in $s^{'}$. The repair operator then re-inserts parts of elements to $s^{'}$ so that $s^{'}$ becomes a complete solution. Let $f(s), f(s^{'}), f(s^{*})$ be the objective of the current solution, the newly obtained solution, and the best-known solution, respectively. At each iteration, an \textbf{acceptance criterion} is used to determine whether the newly obtained solution $s^{'}$ is accepted as the current solution $s$. When the \textbf{termination criterion} is met, the algorithm outputs the optimal solution $s^{*}$. The framework of ALNS is shown in Algorithm \ref{alg_ALNS}. 
 
\begin{algorithm}[H]
	\hspace*{0.02in} {\bf Input:}
	An initial solution $s_{0}$ (Section \ref{sec_IS}), set of destroy operators $\Omega^{-}$, set of  repair operators $\Omega^{+}$.
	\\
	\hspace*{0.02in} {\bf Output:} 
   Best-known solution $s^{*}$.
	\begin{algorithmic}[1]
            \State Initialize $s^{*} \gets s_{0}$, $s^{'} \gets s_{0}$
            \While{Termination Criterion (Section \ref{subsecTC}) not met}
                \State Select a destroy operator (Section \ref{sec_DO}) and a repair operator (Section \ref{sec_RO}) from $\Omega^{-}$ and $\Omega^{+}$ according to the adaptive mechanism in Section \ref{subsecAM}.
                \State $s^{'} \gets repair(destroy(s))$
                \If{$s^{'}$ is better than $s^{*}$}
                    \State $s^{*} \gets s^{'}$
                    \State $s \gets s^{'}$
                \EndIf
                \If{accept($s^{'}$) (Section \ref{subsecAC})}
                    \State $s \gets s^{'}$
                \EndIf
                \State Update the parameters of the adaptive mechanism
            \EndWhile
        \end{algorithmic}
	\caption{(The framework of ALNS)}
	\label{alg_ALNS}
\end{algorithm}

The rest of this section instantiates each step of the ALNS for solving the OTOP-RV.

\subsection{Initial Solution} \label{sec_IS}
	We generate an initial solution for the OTOP-RV via a straightforward heuristic, by sequentially solving single-agent open orienteering problems, once for each agent. This requires the following definition.

\begin{definition}{\textbf{(Marginal gain)}}
	Given a POI $i \in I$ and the number of visits $q_{i}$ already made by the routes, the marginal gain refers to the additional reward obtained from one more visit, namely $\eta_{i}=w_{i}\big((q_{i}+1)^\beta - (q_{i})^\beta\big)$.
\label{def_MG}
\end{definition}
	
	The idea is to sequentially generate a route for each agent by maximizing the marginal gain, which is defined by all previously generated routes. The pseudo-code for the generation of the initial solution is shown in Algorithm \ref{alg_IS}.

\begin{algorithm}[H]
	\hspace*{0.02in} {\bf Input:}
	Number of agents $n_{K}$, the length of time or distance budget $\Delta$, marginal gains $\{\eta_{i}, i \in I\}$
	\\
	\hspace*{0.02in} {\bf Output:} 
   Initial routing solutions $R_{0}$
	\begin{algorithmic}[1]
		\For{$v = 1,...,n_{K}$} 
			\State Generate a new route by solving the Model \eqref{CG_SP_1}-\eqref{CG_SP_10}, which can maximize the marginal gain.
			\State Update the marginal gain of each POI based on the new route.
		\EndFor
   \end{algorithmic}
     
	\caption{(Generation of initial solution)}
	\label{alg_IS}
\end{algorithm}

The single-agent routing problem can be described as a classic orienteering problem. Let $x_{i,j}$ be a binary variable taking the value 1 if the new route directly connects vertex $i$ to vertex $j$, $y_{i}$ be a binary variable taking the value 1 if POI $i$ is visited by the new route. Let $z_{i}$ be the order in which POI $i$ is visited in the new route, and we set $z_{0}=0$ and $z_{n+1}=n+1$ (where $0$ and $n+1$ are deemed the virtual starting and ending depots).

	\begin{equation}
        \label{CG_SP_1}
        \max_{x_{i,j}} \sum_{i \in I}\eta_{i}y_{i}
	\end{equation}

    \begin{eqnarray}
        \label{CG_SP_2}
        \sum_{i \in I}x_{0,i}=1 
        \\
        \label{CG_SP_3}
        \sum_{i \in I}x_{i,n+1}=1
        \\
        \label{CG_SP_4}
        \sum_{(i,l) \in A}x_{i,l} = \sum_{(l,j) \in A}x_{l,j} & & \forall l \in I
        \\
        \label{CG_SP_5}
        \sum_{(i,j) \in A}x_{i,j} = y_{i}  & & \forall i \in I
        \\
        \label{CG_SP_6}
        \sum_{(i,j) \in A}t_{ij}x_{i,j} \leq \Delta
        \\
        \label{CG_SP_7}
         z_{i}-z_{j}+1 \leq n(1-x_{i, j}) & & \forall (i,j) \in A
        \\
        \label{CG_SP_8}
         1 \leq z_{i} \leq n & & \forall i \in I
        \\
        \label{CG_SP_9}
        x_{i,j} \in \{0,1\}  & & \forall (i,j) \in A
        \\
        \label{CG_SP_10}
        y_{i} \in \{0, 1\}  & & \forall i \in I
    \end{eqnarray}
	
	The objective function \eqref{CG_SP_1} maximizes the total marginal gain collected by the new route. Constraints \eqref{CG_SP_2} and \eqref{CG_SP_2} ensure that the route starts from node $0$ and ends with node $n+1$ (i.e. the depots). Constraints \eqref{CG_SP_4} expresses flow conservation. Constraints \eqref{CG_SP_5} relates variables $x_{i,j}$ to the indicator variable $y_{i}$. Constraint \eqref{CG_SP_6} ensures that the time or distance budget is not exceeded. Constraints \eqref{CG_SP_7} and \eqref{CG_SP_8} act as sub-tour elimination constraints \citep{MTZ1960}.

	 We note that the Model \eqref{CG_SP_1}-\eqref{CG_SP_10} is an orienteering problem, which already has mature solutions. In this study, we solve this model by the heuristic algorithm proposed by \cite{Tsiligirides1984}.

\subsection{Destroy Operators}\label{sec_DO}
	This section describes four removal heuristics: random removal, worst removal, related removal, and route removal. All four heuristics take the current solution $s$ as input. The output of the heuristic is a temporary solution $s_{temp}^{'}$ after applying the destroy operator, which removes some points from the current solution $s$.

\subsubsection{Random removal}
	In the random removal heuristic, we randomly remove $m\%$ of the points (any decimals will be rounded to the nearest integer) from the route of each agent using a uniform probability distribution.

\subsubsection{Worst removal}\label{sec_WR}
	As proposed in \cite{RP2006}, the worst removal heuristic removes points to obtain the most savings. For the OTOP-RV, we need to consider the trade-off between collecting rewards and the time/distance required. Specifically, we aim to remove the POIs that have limited impact on the total reward, while resulting in the most time/distance savings.

	Let $\mathcal{V}$ be the set of visited POIs in the current solution $s$. Given a POI $i\in\mathcal{V}$, we use $f(s)$ to represent the reward of the current solution $s$, and $T(s)$ to represent the total time/distance used by all agents in the current solution $s$. Next, we define the reward lost by removing $i$ from the current solution $s$ as $r(i,s)=f(s)-f_{-i}(s)$, and the time/distance saved as $t(i,s)=T(s)-T_{-i}(s)$. We define the value of each POI in the current solution $s$ as $value(i,s)$:
	\begin{equation}
		\label{PV}
		value(i,s) = \frac{r(i,s)}{t(i,s)}
	\end{equation}

\noindent Building on such a notion, this operator iteratively removes POIs with relatively low values. Following \cite{RP2006}, we introduce random factors to avoid situations where the same points are removed over and over again. In addition, randomization is applied in a way that favors the selection of POIs with lower values. The pseudo-code for the worst removal heuristic is shown in Algorithm \ref{alg_WR}.

\begin{algorithm}[H]
	\hspace*{0.02in} {\bf Input:}
	A current solution $s$, proportion $m$\% of point removal, parameter $p >0$
	\\
	\hspace*{0.02in} {\bf Output:} 
   Temporary solution $s_{temp}^{'}$ after destroy operator 
	\begin{algorithmic}[1]
		\State Calculate the number of POIs to be removed based on $m$\%, denoted as $n_{p}$;
		\For{$n = 1,...,n_{p}$} 
			\State List of all visited POIs in the current solution $s$, denoted as $\mathscr{V}$;
			\State $\mathscr{V}_{new}$ = Sorted list of $\mathscr{V}$ in an increasing order of $value(i,s)$, $i\in\mathcal{V}$;
			\State Draw a random number $y$ from the uniform distribution $U(0, 1)$;
			\State Let $x$ be the $\lfloor y^{p}|\mathscr{V}| \rfloor$-th element in the ordered set $\mathscr{V}_{new}$ ($\lfloor\cdot \rfloor$ is the floor operator);
			\State Remove point $x$ from the current solution $s$.
		\EndFor
   \end{algorithmic}
	\caption{(Worst removal)}
	\label{alg_WR}
\end{algorithm}

\subsubsection{Related removal} 
	The purpose of the related removal heuristic is to remove a set of points that, in some way, are closely related and hence easy to interchange among different routes during repairs \citep{Pisinger2007}. For the OTOP-RV, we remove some adjacent POIs, with the understanding that points closer to each other are more likely to be interchanged. Specifically, we first randomly select a POI $o \in I$ as the center point, and then remove the $m\%$ points closest to $o$ in the current solution $s$.

\subsubsection{Route removal}
     This operator is commonly used in the vehicle routing problem \citep{Demir2012}. For the OTOP-RV, the operator randomly removes $\left[n_K\times m\%\right]$ of the routes from the solution $s$, where $n_K$ is the number of routes ($\left[\cdot \right]$ is the rounding operator).
     
\subsection{Repair Operators} \label{sec_RO}
    This section describes some insertion heuristics: greedy insertion and (k)-regret insertion. Insertion heuristics are typically divided into two categories: sequential and parallel. The difference between the two is that the former builds one route at a time while the latter constructs several routes simultaneously \citep{Potvin1993}. The insertion heuristics adopted in this paper are all parallel, which will be used to repair the temporary solution $s_{temp}^{'}$ following the destroy operator.  The output of the insertion heuristic is a new feasible solution $s^{'}$.

\subsubsection{Greedy insertion}\label{subsubsecgi}

    \begin{definition}{\textbf{(Minimum-cost position)}} \label{D_MCP}
    	Given a POI $i \in I$ and a route $R_{k}$, there are $|R_{k}|+1$ possible insertion positions for this POI. Let $\Delta t_{ik}(m)$ be the extra time/distance incurred by inserting $i$ into the $m$-th position of $R_k$. We define the minimum-cost position to be $m_{ik}^{*}=\mathop{\arg\min}\limits_{m=\{1,..., |R_{k}|+1\}}\Delta t_{ik}(m)$ and the minimum cost as $\Delta t_{ik}(m_{ik}^{*})$.
    \end{definition}

\noindent Let $\Delta r_{i, k}$ be the increased reward after inserting point $i$ into the route of agent $k$ at the minimum-cost position. If the POI $i$ cannot insert to the route of agent $k$, we set $\Delta r_{i, k} = 0$

    \begin{definition}\label{D_VP}
        Similar to the reward efficiency defined in \eqref{RE}, the reward efficiency of inserting point $i$ is defined as the extra reward collected per unit time/distance required to accommodate such insertion: 
$$
\Psi_{i,k}\doteq {\Delta r_{i,k}\over \Delta t_{ik}(m_{ik}^*)}
$$
\noindent Then, we define the maximum insertion value of $i$ to be
        \begin{equation}
		\label{PV_2}
		v(i)\doteq\max_{k \in K}\Psi_{i,k}
	\end{equation}
    \end{definition}

The greedy insertion operator iteratively selects and inserts a candidate POI $i$ that has the maximum insertion value $v(i)$. This process continues until no more POIs can be inserted into any route.

\subsubsection{(k)-Regret insertion}
    The regret heuristic tries to improve upon the greedy heuristic by incorporating look-ahead information when selecting the point for insertion \citep{RP2006}. In prose, the algorithm performs the insertion that will be most regretted if it is not done now. Recalling the notions from Section \ref{subsubsecgi}:
    \begin{eqnarray}
        \label{kRI_1}
        \Delta v_{ik}(m) = \frac{\Delta r_{i, k}}{\Delta t_{ik}(m)} \qquad \forall i \in I, \quad k \in K,\quad m = \{1, ..., |R_{k}|+1\}
    \end{eqnarray}

\noindent We sort the list $\big\{\Delta v_{ik}(m):  i \in I,  k \in K, 1\leq m \leq |R_{k}|+1\big\}$ in descending order to obtain $\big\{\Delta v^{(n)}_{ik}(m): i \in I,  k \in K, 1\leq m \leq |R_{k}|+1\}\big\}$ where $n$ indexes the ordered  elements in the list. Then, the $k$-regret heuristic chooses to insert the POI that maximizes $c_{k}^{*}$ in each step.

    \begin{equation}
        \label{kRI_2}
        c_{k}^{*}= \sum_{j=1}^{k}(\Delta v_{ik}^{j}(m) - \Delta v_{ik}^{1}(m))
    \end{equation}

\noindent
where $c_{k}^{*}$ is the $k$-regret value. This process continues until no POI can be inserted into any route.

\subsection{Adaptive Mechanism}\label{subsecAM}

	We defined, in Section \ref{sec_DO},  four destroy operators (random, worst, related, and route removal), and in Section \ref{sec_RO} a class of repair operators (greedy insertion, $k$-regret). This section explains a strategy to adaptively select the destroy and repair operators at each iteration of ALNS.
	
	We assign weights to different operators and implement the classical roulette wheel mechanism \citep{RP2006}. If we have $L$ operations with weights $\mu_{l}$, $l \in \{1,2,...,L\}$, we select operator $l$ with probability $p_{l}$ given as:
\begin{equation}
        \label{Choose_operators}
        p_{l}= \frac{\mu_{l}}{\sum_{l=1}^{L}\mu_{l}}
\end{equation}

	Initially, we set the same weight for each operator. Then, we continuously update the weight by keeping track of a score for each operator, which measures how well the heuristic has been performing recently. Operators with higher scores have a higher probability of being selected during the iteration. The entire search is divided into several segments, indexed by $h$. A segment contains a few iterations of the ALNS, hereafter defined to be $\delta$ iterations. The score of all operators is set to zero at the start of each segment. The scores of the operators are updated according to the following rules:
\begin{itemize}
\item[(a)] if the score of the new solution is better than the best-known solution, the operator score increases by $\sigma_{1}$; 
\item[(b)] if the score of the new solution is worse than the best-known solution but better than the current solution, the operator score increases by $\sigma_{2}$; 
\item[(c)] if the score of the new solution is worse than the current solution, but is accepted, the operator score increases by $\sigma_{3}$; 
\item[(d)] if the new solution is not accepted, the operator score increases by $\sigma_{4}$. 
\end{itemize}

It is reasonable to set $\sigma_{1} > \sigma_{2} > \sigma_{3} > \sigma_{4}$. At the end of each segment, we calculate new weights using the recorded scores. Let $\mu_{l,h}$ be the weight of operator $l$ used in segment $h$, and the probability $p_{l,h}$ is calculated according to \eqref{Choose_operators}. At the end of each segment $h$ we update the weight for operator $l$ to be used in segment $h+1$ as follows:

\begin{equation}
        \label{Update_weights}
        \mu_{l,h+1}=  \mu_{l,h}(1-\varepsilon) + \varepsilon \frac{\pi_{l,h}}{\theta_{l,h}}
\end{equation}

\noindent
where $\pi_{l,h}$ is the score of operator $l$ obtained during the last segment and $\theta_{l,h}$ is the number of times operator $l$ was used during the last segment. The reaction factor $\varepsilon$ controls how quickly the weight adjustment algorithm reacts to changes in the effectiveness of the operators.

\subsection{Acceptance Criterion}\label{subsecAC}

	A simulated annealing framework is used to decide whether to accept the newly obtained solution $s^{'}$ given the current solution $s$. If the objective value $f(s^{'}) > f(s)$, then $s^{'}$ is accepted as the current solution. Otherwise, the probability of $s^{'}$ being accepted is $P(s \leftarrow s^{'})$.

\begin{equation}
        \label{SA}
        P(s \leftarrow s^{'})= \exp\left\{\frac{f(s^{'})-f(s)}{T}\right\}
\end{equation}

\noindent
where $T >0$ is the temperature, which begins with $T_{\text{max}}$ and  decreases at every iteration according to $T = T \cdot c$, where $c\in (0,1)$ is the cooling rate. When $T$ is below $t_{\text{min}}$, take reheating measures and set $T$ to $T_{\text{max}}$.

\subsection{Termination Criterion}\label{subsecTC}

	The algorithm stops after (1) a maximum number of iterations $N_{1}^{\text{max}}$; or (2) a certain number of non-improving iterations $N_{2}^{\text{max}}$.

\section{Discussion and extension}\label{secDE}
This section provides some discussion on the OTOP-RV when modeling various real-world problems. 

\subsection{Completeness of graphs}
The OTOP-RV proposed in this paper takes as input the POI set $I$, their weights in the objective $\{w_i:\,i\in I\}$, and the adjacency matrix. In practice, the adjacency matrix is replaced by the travel time (or distance) matrix $\{t_{ij}:\, i, j\in I\}$, where $t_{ij}$ denotes the travel time (or distance) from vertex $i$ to $j$. 

\begin{itemize}
\item In application scenarios like UAV routing, the spatial domain of the problem is a subset of a Euclidean space, in which any two vertices are directly connected. In this case, any element of the adjacency matrix $t_{ij}$ is a finite number. In other words, all the POIs can be seen as vertices in a complete graph. 

\item In application scenarios like car routing, the spatial domain is a road network represented as a directed graph $G(\mathcal{V},\,\mathcal{A})$, where $\mathcal{V}$ is the set of vertices and $\mathcal{A}$ is the set of arcs, we set $t_{ij}=\infty$ whenever the arc $(i, j)\notin\mathcal{A}$. In other words, road networks can be seen as incomplete graphs.
\end{itemize}

\subsection{Arcs as POIs}	
In some applications based on road networks, the rewards are collected on arcs instead of nodes. Examples include road surface sprinkling (aiming at reducing fugitive dust), or road roughness monitoring (for regular maintenance and repair). Such tasks need to be repeatedly performed, and the weights of these arcs are heterogeneous. These problems fall within the purview of OTOP-RV because one can simply augment the original network with artificial nodes (with rewards assigned) in the middle of relevant arcs. For example, if arc $(i, j)\in\mathcal{A}$ carries rewards, we insert a node $k$ to form two new arcs $(i,k)$ and $(k, j)$, and set $t_{ik}=t_{kj}={1\over 2} t_{ij}$. Then, the weight and reward of arc $(i,j)$ in the original network are transferred to node $k$ in the augmented network. Such a simple technique converts arc-based OTOP-RV to node-based OTOP-RV, which can be solved with the proposed algorithms. 

\section{Computational studies} \label{sec_computational_study}
	All the computational performances reported below are based on a Microsoft Windows 10 platform with Intel Core i9 - 3.60GHz and 16 GB RAM, using Python 3.8 and Gurobi 9.1.2.

\subsection{Instances sets and parameter settings}
\begin{figure}[h]
       \centering
       \includegraphics[width=0.85\textwidth]{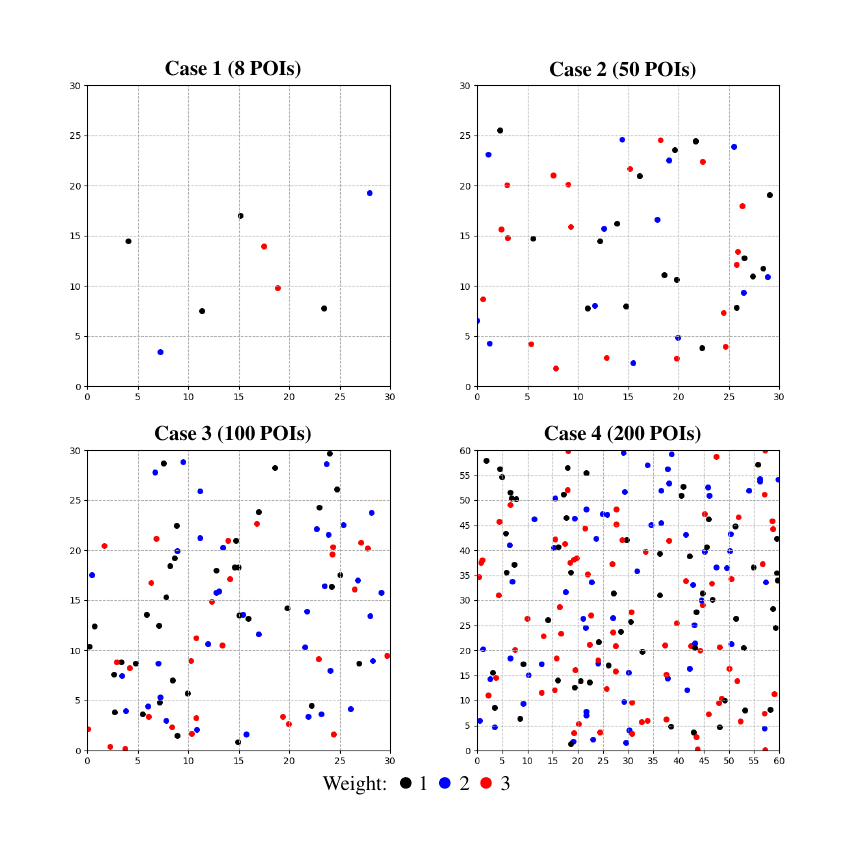}
       \caption{Four test datasets with randomly generated locations and weights of POIs.}
       \label{fig_four_case}
   \end{figure}

	As the proposed OTOP-RV is a new variant of TOP, no existing benchmark instances
are available to evaluate the solution’s performance. Therefore, we randomly generate 20
instances containing 8, 50, 100, and 200 POIs, as shown in Figure \ref{fig_four_case}. Each set of problems consists of 5 instances. All the results reported below for each case are averaged over 5 instances. Note that in this case, the travel times $t_{ij}$ are directly calculated as the Euclidean distance between two points. In other words, those POIs depicted in Figure \ref{fig_four_case} are vertices in complete graphs. 

	The performance of the ALNS algorithm is sensitive to the settings of parameters. In
our tests, the range and interval of each parameter are determined by extensive tests. In
each test, we draw an initial value from the range of each parameter and tune it gradually.
The final settings of the parameters are as follows. The destroy operators are parameterized by $m$, we set $m=0.4$. The insertion heuristics in the repair operators are parameter-free. The weight adjustment algorithm is parameterized by $\sigma_{1}$, $\sigma_{2}$, $\sigma_{3}$, $\sigma_{4}$ and learning rate $\varepsilon$. We set $(\sigma_{1}$, $\sigma_{2}$, $\sigma_{3}$, $\sigma_{4}, \varepsilon)=(20, 10, 3, 0, 0.7)$. To manage the acceptance criterion we use three parameters, $T_{\text{max}}$, $t_{\text{min}}$ and $c$. The start temperature $T_{\text{max}}$ is determined based on the value of the initial solution $f(s_{0})$ and the formula is $T_{\text{max}} = 0.05\times \frac{f(s_{0})}{\ln(2)}$, which allows for a 50\% probability that non-improving solutions generated at the initial solution temperature is accepted. In addition, we set $t_{\text{min}}=0.1$ and $c=0.95$. We also need to determine $N_{1}^{\text{max}}$ and $N_{2}^{\text{max}}$ that govern algorithm termination. We set $N_{1}^{\text{max}}$ to a large number 2000, and mainly rely on $N_{2}^{\text{max}}$ to terminate the algorithm. The value of parameter $N_{2}^{\text{max}}$ set to 200.

\subsection{Algorithm performance}\label{subsecAP}
	Because we are the first to explore the OTOP-RV, there are no existing solution approaches to benchmark with. so we use the following three benchmarks for comparison.

\begin{itemize}
	\item[(1)] \textbf{Gurobi solver:} For small cases, we prove the effectiveness of the ALNS algorithm by comparing it with- the exact solution of the solver.
	\item[(2)] \textbf{Greedy algorithm:} For large cases, we compare the ALNS algorithm with a greedy algorithm. The idea is to sequentially select the next POI that maximizes the reward efficiency among all feasible POIs. For a detailed description of this method, see the Appendix, and the pseudo-code for the greedy algorithm is shown in Algorithm \ref{alg_greedy}.
	\item[(3)] \textbf{Sequential OP:} This method aims to sequentially generate a route for each agent by maximizing the marginal gain, which is defined by all the previously generated routes. The pseudo-code is shown in Algorithm \ref{alg_IS}.

\end{itemize}

We evaluated the effectiveness of the ALNS algorithm through testing on Case 1 (8 POIs) as shown in Figure \ref{fig_four_case}. The test results under different parameters are shown in Table \ref{tab_small_case}. As the $n_{K}$ and $\Delta$ increases, the computational time of Gurobi grows faster than that of ALNS. For small cases, the ALNS algorithm can find the same optimal solutions as Gurobi but in a shorter time.

	\setlength\LTleft{0pt}
	\setlength\LTright{0pt}
	\begin{longtable}[c]{cc|ccc|ccc}
		\caption{Computational results comparing the Gurobi and ALNS in a small case. The computational time limit of Gurobi is set to 7200s.}
            \label{tab_small_case} 
            \\
	    \hline
		\multirow{2}{*}{$n_K$}	&	\multirow{2}{*}{$\Delta$} & \multicolumn{3}{c|}{Gurobi}& \multicolumn{3}{c}{ALNS}  \\
		\cline{3-8} 
		& & obj & gap & CPU(s) &  obj  & CPU(s)  & Increase
		\\
		\hline
		2 & 20 & 13.5  & 0.0\%  & 21  & 13.5  &  1.0   & 0.0\%    
		\\
		\hline
		2 & 30 & 16.0  & 0.0\%  & 75   & 16.0  &  1.2   & 0.0\%    
		\\
		\hline
		2 & 40 & 17.9  & 0.0\%  & 220   & 17.9  &  2.0  & 0.0\%    
		\\
		\hline
		3 & 20 & 17.2  & 0.0\%  & 477  & 17.2  &  2.0   & 0.0\%    
		\\
		\hline
		3 & 30 & 20.8  & 0.7\%  & 3004   & 20.5  &  3.2   & -1.4\%    
		\\
		\hline
		3 & 40 & 22.8  & 1.1\%  & 6177   & 22.8  &  4.4   & 0.0\%    
		\\
		\hline
		4 & 20 &  20.1 & 23.4\%  &   7200&  20.1 &   5 & 0.0\%    
		\\
		\hline
		4 & 30 &  23.7 & 18.9\%  &  7200 & 23.8  &   5.4 & 0.4\%    
		\\
		\hline
		4 & 40 & 26.2  & 7.7\%  & 7200 & 26.2  &  6.8   & 0.0\%    
		\\
		\hline
	\end{longtable}

\begin{figure}[h]
       \centering
       \includegraphics[width=\textwidth]{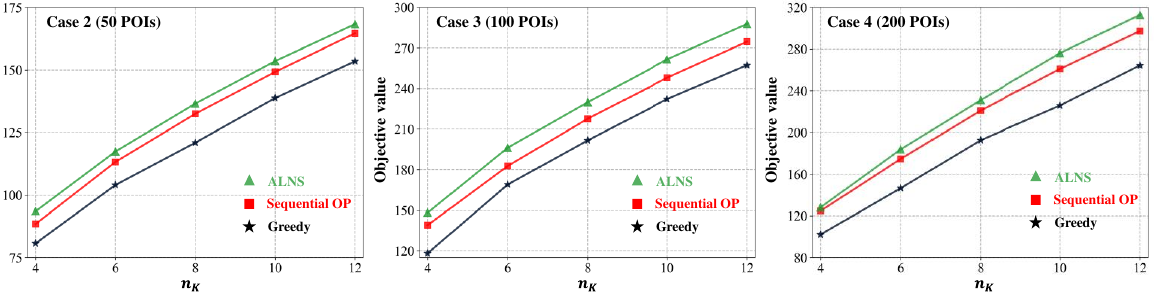}
       \caption{Objective values in three test cases with different values of $n_K$.}
       \label{fig_benchmark}
   \end{figure}

We compare the objective value of the greedy algorithm, sequential OP, and ALNS on the three datasets (Case 2, 3, 4) as shown in Figure \ref{fig_four_case}. The test results under different $n_{K}$ (number of routes) are shown in Figure \ref{fig_benchmark} and Table \ref{tab_benchmark}. It is shown that: 
\begin{itemize}
\item[(1)] The ALNS significantly outperforms the greedy algorithm in objective value, with improvements ranging from 9.7\% to 25.4\%; 
\item[(2)] The ALNS can further improve the solution of sequential OP by neighborhood search, with improvements ranging about 6\%;
\item[(3)] Such improvements are less pronounced for larger $n_K$ (number of routes), because of (i) the diminishing marginal gain by design; and (ii) the relatively high saturation of agent routes.
\end{itemize}

\setlength\LTleft{0pt}
	\setlength\LTright{0pt}
	\begin{longtable}[c]{cc | c | cc | ccc}

	\caption{Comparison of objective values among greedy algorithm, sequential OP and ALNS. The `Increase' column represents improvements over the greedy algorithm.}
        \label{tab_benchmark} 
        \\
	    \cline{2-7}
         		&\multirow{2}{*}{$n_K$} & \multicolumn{1}{c|}{Greedy} & \multicolumn{2}{c|}{Sequential OP} & \multicolumn{2}{c}{ALNS}  \\
        \cline{3-7} 
         & & Obj. & Obj. & Increase & Obj. & Increase  
        \\
        \hline
        \multirow{4}{*}{{\bf Case 2}} & 4   & 80.7  & 88.4  & 9.6\%  & 93.5 &  15.9\%    
        \\
        \multirow{4}{*}{50 POIs}  & 6   & 104.0  & 113.2  & 8.9\%  & 117.4 & 12.9\%    
        \\
       \multirow{4}{*}{$\Delta=30$}   & 8   & 120.9  & 132.5   & 9.5\%  & 136.5 &  12.9\%    
        \\
          & 10   & 138.7  & 149.2 & 7.6\%  & 153.6 &   10.7\%    
        \\
          & 12  & 153.5  & 164.6  & 7.3\%  & 168.3 &   9.7\%    
        \\
        \hline
        \multirow{4}{*}{{\bf Case 3}}  & 4   & 118.2  & 139.0   & 17.6\%  & 148.2 &  25.4\%    
        \\
          \multirow{4}{*}{100 POIs} & 6   & 169.0  & 182.7  & 8.1\%  & 196.2 &   16.1\%    
        \\
        \multirow{4}{*}{$\Delta=30$}  & 8   & 201.7  & 217.7  & 7.9\%  & 229.9 &   14.0\%    
        \\
         & 10   & 232.2  & 248.0   & 6.8\%  & 261.6 &  12.7\%    
        \\
         & 12  & 257.3  & 274.6 & 6.7\%  & 287.7 &  11.8\%        
        \\
        \hline
        \multirow{4}{*}{{\bf Case 4}} & 4   & 102.2  & 124.8  & 22.1\%  & 128.2 &   25.4\%    
        \\
       \multirow{4}{*}{200 POIs}  & 6   & 146.6 & 174.7 & 19.2\%  & 183.8 &  25.4\%    
        \\
       \multirow{4}{*}{$\Delta=30$}  & 8   & 192.6  & 221.0   & 14.7\%  & 231.0 &   20.0\%    
        \\
          & 10   & 225.8 & 261.0  & 15.6\%  & 276.0 &  22.2\%    
        \\
          & 12   & 264.5  & 297.3  & 12.4\%  & 312.8 &  18.3\%    
        \\
        \hline
    \end{longtable}

\begin{figure}[H]
       \centering
       \includegraphics[width=\textwidth]{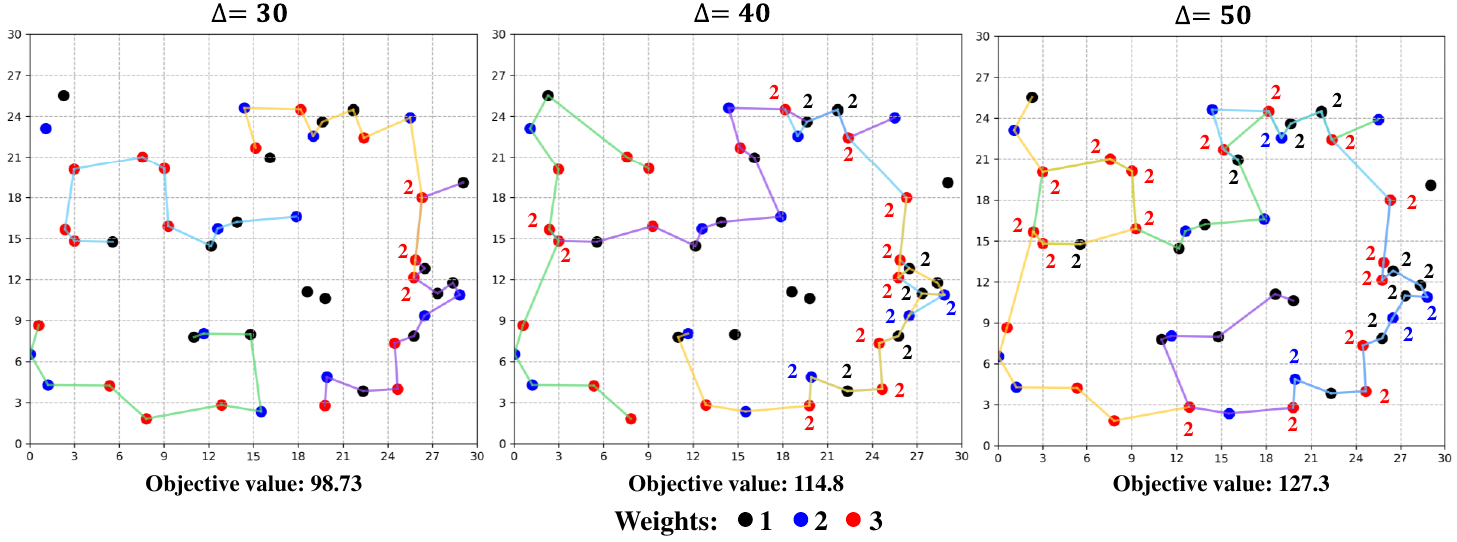}
       \caption{Solution visualization of Case 2 (50 POIs) with four routes. $\Delta$ is the distance budget for each route. Nodes visited more than once are labeled with the number of visits.}
       \label{fig_solution1}
   \end{figure}

    Figure \ref{fig_solution1} provides a visualization of the OTOP-RV solution for Case 2, with four routes and different distance budget $\Delta$. It can be seen that, even for a small-scale problem, the solution displays some complexity, especially for larger $\Delta$ where the routes have considerable overlap. Moreover, such overlap took place at high-value nodes, which reflects the effectiveness of our proposed algorithms. Such solutions cannot be obtained via conventional TOPs or their variants.
    
\subsection{Real-world case study}\label{subsecRWCS}
	The OTOP-RV is demonstrated in a real-world case of Volatile Organic Compounds (VOCs) monitoring in the Longquanyi District, Chengdu, China. VOCs include a wide variety of chemicals, some of which have adverse health effects and act as catalysts of processes that form PM$_{2.5}$ and O$_3$ \citep{Mishra2023, Yang2020}. VOCs are emitted from the manufacturing activities of 141 factories in Longquanyi, and a single VOCs sensing vehicle, operated by a local environmental protection agency, needs to perform regular monitoring tasks of the area by moving along designated routes. 
	
\begin{figure}[h]
       \centering
       \includegraphics[width=1.0\textwidth]{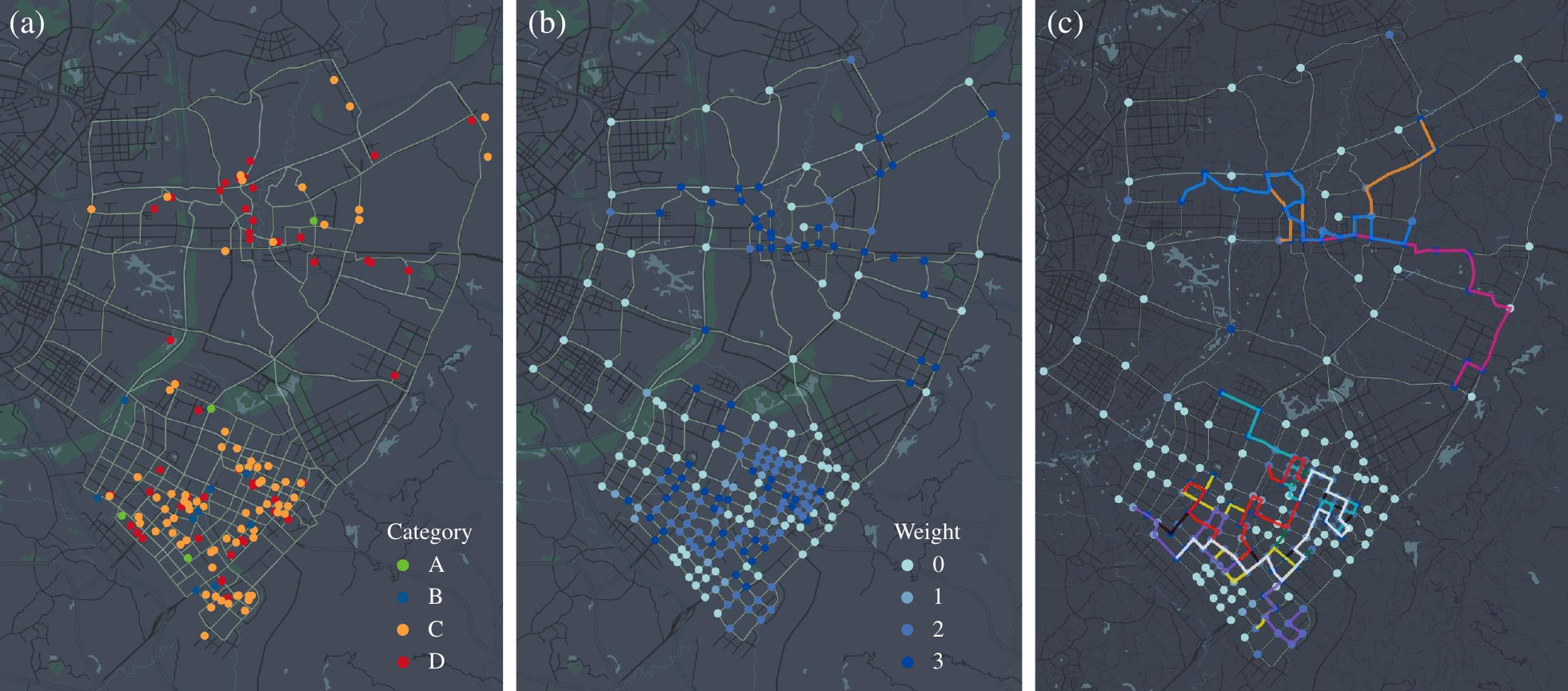}
       \caption{(a): Road network and distribution of factories of Longquanyi District, Chengdu. (b): The sensing weights of all the road network nodes. (c) Visualization of the 10 routes generated by the ALNS algorithm.}
       \label{fig_factory}
   \end{figure}

	The 141 factories are categorized into four classes: A, B, C, and D, where class A has the lowest sensing priority and D has the highest; see Figure \ref{fig_factory}(a). Based on the locations of these factories and their relative positions to the road network, their sensing priorities are transformed into the weights of the network nodes, which are treated as the POIs in our model; see Figure \ref{fig_factory}(b).

	To assess the validity of the solution in a mobile sensing context, we assume that the sensing vehicle needs to perform two routes per day, each within a distance budget of 20 km\footnote{These are in line with real-world operations in Longquanyi}. Therefore, 10 routes need to be generated for a week's (5 working days) monitoring task. This is formulated as an OTOP-RV and solved by the ALNS. Figure \ref{fig_factory}(c) shows the 10 routes generated by the algorithm. In this solution, 7 routes are located in the southern part of the area where the majority of the POIs are located. Moreover, Table \ref{tab_cover_times} shows that many of the high-value POIs (with higher weights) are covered by multiple routes, which means those factories with higher sensing priority are more frequently visited, which is a desired feature of the routing plan.

\begin{table}[h]
\centering
\caption{Number of POIs meeting different coverage times}
\label{tab_cover_times}
\resizebox{\columnwidth}{!}{%
\begin{tabular}{|c|c|c|c|c|c|c|}
\hline
\multirow{2}{*}{Weight} & \multicolumn{6}{c|}{Number of POIs}                                                                                                  \\ \cline{2-7} 
 & Covered by 4 routes & Covered by 3 routes & Covered by 2 routes & Covered by 1 route & Not covered & Sum \\ \hline
3    & 6  & 18 & 17 & 19 & 2  & 62 \\ \hline
2  & 13 & 20 & 16 & 13 & 4  & 66 \\ \hline
1    & 0  & 1  & 4  & 4 & 3  & 12 \\ \hline
0    & 0 & 0 & 2  & 4  & 76 & 82 \\ \hline
\end{tabular}%
}
\end{table}

\section{Conclusion} \label{sec_conclusion}

This paper proposes an extension of the team orienteering problem (TOP), named OTOP-RV, by allowing a POI to be visited multiple times to accumulate rewards. Such an extension is important because many real-world applications such as mobile sensing and disaster relief require that the POIs are repeatable visited, and such needs are heterogeneous. To model such scenarios, we present the OTOP-RV as a nonlinear integer program and propose an ALNS-based heuristic algorithm. The following findings are made from extensive numerical tests.
 
\begin{itemize}
\item The ALNS is effective in finding good-quality solutions, where high-value POIs are visited more frequently. This is not achievable through solution approaches for conventional TOPs. 

\item For small cases, the ALNS algorithm can find the same optimal solutions as Gurobi but in a shorter time. 

\item For large cases, the ALNS significantly outperforms the greedy algorithm in objective value, with improvements ranging from 9.7\% to 25.4\%.

\item The ALNS can further improve the solution of sequential OP by neighborhood search, with improvements ranging about 6\%. 
 
\item The OTOP-RV is transferrable to treat scenarios where the demands are concentrated on arcs instead of vertices/nodes. The resulting model can be applied to anti-dust (road sprinkler) operations or road surface condition monitoring.
\end{itemize}

\section*{Appendix}
\subsection*{A Greedy algorithm for the OTOP-RV}
This greedy algorithm generates the routes of each agent in sequence. The calculation of the marginal gain $\eta_{i}$ of point $i$ as shown in Definition \ref{def_MG}. We begin with the POI $j \in I$ with the highest marginal gain as the starting location for the agent. The idea is to sequentially select the next POI $i^{*}$ that maximizes the reward efficiency $\psi_{ji}$ among all feasible POIs:

$$
		\label{RE}
		i^{*} = \mathop{\arg\max}\limits_{i}\psi_{ji}= \mathop{\arg\max}\limits_{i}\frac{\eta_{i}}{t_{ji}}
$$
                        
\noindent In prose, the reward efficiency refers to the reward collected per unit spending of time (or distance). The pseudo-code for the generation of the initial solution is shown in Algorithm \ref{alg_greedy}.

\begin{algorithm}[H]
	\hspace*{0.02in} {\bf Input:}
	Number of agents $n_{K}$, the length of time or distance budget $\Delta$, marginal gains $\{\eta_{i}, i \in I\}$
	\\
	\hspace*{0.02in} {\bf Output:} 
   Initial routing solutions $R_{0}$
	\begin{algorithmic}[1]
		\For{$v = 1,...,n_{K}$} 
			\State Initialize $t = 0$, list of candidate points $I_{temp}=I$;
			\State Set the starting point of agent $v$ as $s_{v} = \mathop{\arg\max}\limits_{i} \eta_{i}$;
			\While {$t < \Delta$}
				\State Find the next visiting point  $i^{*}=\mathop{\arg\max}\limits_{i}\psi{(i)}$;
				\State Add the point $i^{*}$ to the agent route $R_{0}^{v}$ and update the marginal gains $\{\eta_{i}| i \in I \}$;
				\State Update $t = t + t_{j, i^{*}}$;
				\State Set the current point is $j \leftarrow i^{*}$;
				\State Remove $j$ from the list of candidate points $I_{temp}$.
			\EndWhile
		\EndFor
   \end{algorithmic}
     
	\caption{(Greedy heuristic algorithm)}
	\label{alg_greedy}
\end{algorithm}

\section*{Data availability}
The datasets that support the findings of this study are openly available in \url{https://github.com/Shenglin807/OTOP-RV-Dataset}.

\section*{Acknowledgement}
This work is supported by the National Natural Science Foundation of China through grants 72071163 and 72101215, and the Natural Science Foundation of Sichuan Province through grants 2022NSFSC0474 and 2022NSFSC1906.

\end{document}